\documentclass[11pt]{article}

\usepackage{amssymb,mathtools}

 \usepackage{amsthm}
 \usepackage{dsfont}
 
\newtheorem{theorem}{Theorem}[section]

\newtheorem{lemma}{Lemma}[section]
\newtheorem{definition}{Definition}[section]
\newtheorem{conjecture}{Conjecture}[section]
\newtheorem{example}{Example}[section]
\usepackage{lineno}
\newcommand{\R}{\mathds{R}}

\title{Two curve Chebyshev approximation and its application to signal clustering}

\author{Dr~N. Sukhorukova, e-mail: nsukhorukova@swin.edu.au\\
Faculty of Science, Engineering and Technology,\\
 Swinburne University of Technology, \\
 PO Box~218, Hawthorn, Victoria, \textsc{Australia} and\\ 
 Centre for Informatics and Applied Optimization,\\
  Federation University Australia}

\begin{document}
\maketitle
\begin{abstract}
In this paper we  extend a number of important results of the classical Chebyshev approximation theory to the case of simultaneous approximation of  two or more functions. The need for this extension is application driven, since such kind of problems appears in the area of curve (signal) clustering. In this paper we propose a new efficient algorithm for signal clustering 
and develop a procedure that allows one to reuse the results obtained at the previous iteration without recomputing the cluster centres from scratch. This approach is based on the extension of  the classical de la  Vall{\'e}e-Poussin's procedure originally developed for polynomial approximation. In this paper, we also develop necessary and sufficient optimality conditions for two curve Chebyshev approximation, that is our core tool for curve clustering. These results are based on application of nonsmooth convex analysis. 
\end{abstract}







%
%
%

\section{Introduction}\label{sec:introduction}


In signal processing, there is often a need for constructing cluster signal prototypes. Cluster prototypes can be interpreted as summary curves that may replace the whole group of signal segments (clusters of curves), since it is assumed that the curves from the same cluster are similar to each other in one way or another. Signal prototypes may be used  for characterising the structure of the signal segments and their common features (data analysis) and also for reducing the amount of information to be stored (data compression). 

$k$-means is a very fast method developed for clustering points in $\R^n$. The name ``$k$-means'' was first proposed by James MacQueen in 1967~\cite{kmeansfirst}. This method is based on the minimisation of the total dissimilarity function: the sum of squares of the Euclidean distances between the points and the corresponding cluster centres. The theoretical and computational properties of this method as well as its applicability to data analysis, signal processing and data mining problems have been studied for decades~\cite{spath80, BagUgonCluster}. 

The classical $k$-means method contains two steps. Firstly, we assign each point to the cluster with the nearest centre. Secondly, for each cluster, we recompute the centres by minimising the total sum of dissimilarities within each cluster. Then we repeat these two steps until none of the points changes its clustering membership. In the case of the classical $k$-means algorithm, the second step is reduced to computing the barycentre  of the cluster and therefore the algorithm is very fast.   

The $k$-means method can be easily applied to curve clustering if the curves are discretised. This can be done, for example, by treating each time moment $t_i$, $i=1,\dots,n$ as a separate coordinate in $\R^n$, where $n$ may be large. 

There have been several modifications of this method, among them the $k$-medoid method~\cite{kaufman:clustering1990}, where the dissimilarity function is based on other types of distances. A comprehensive review of $k$-medoid types of methods can be found in~\cite{kmedoid}.  Most $k$-medoid methods are slower than the classical $k$-means. At the same time, they may be more appropriate for some specific models and therefore there is a strong need for studying these methods as well.

It is desirable for a $k$-medoid algorithm, that any cluster prototype is an accurate approximation of each member of the cluster. On the top of this, it is important that the process of recomputing cluster prototypes, when groups of signals move from one cluster to another, is not computationally expensive. 

 In this paper we suggest a uniform (Chebyshev) approximation based model and therefore our method is from the $k$-medoid groups of methods. This is a convex optimisation problem. There are several advantages of the proposed implementation model. First of all, it  provides an accurate approximation to the group of signals. Second, this problem can be reformulated as a linear programming problem, that can be solved efficiently. Finally, the proposed  approach  allows one to compute prototype updates without recomputing from  scratch.   
 
The paper is organised as follows. In section~\ref{sec:math_formulation}, we demonstrate that this problem can be formulated as a linear programming problem and study possible ways for solving this problem efficiently. Then, in section~\ref{sec:opt_cond}, we develop the necessary and sufficient optimality conditions for curve clustering. 
In section~\ref{sec:modelling}, we extend the classical de la  Vall{\'e}e-Poussin's procedure to the case of two curve approximation. This procedure plays a vital role in cluster prototype computations. Finally, in section~\ref{sec:conc}, we comment on the results and underline our future research directions.

\section{Mathematical formulation}\label{sec:math_formulation}
\subsection{Prototype construction}
Assume that there is a group of $l$ signals $S_1(t),\dots,S_l(t)$, whose values are measured at discrete time moments $$t_1,\dots, t_N,~t_i\in[a,b],~i=1,\dots,N.$$
We suggest to construct a group prototype in the form $$S({\bf A},t)=\sum_{i=0}^{n}a_ig_i(t).$$ The functions $\{g_i(t)\}_{i=0}^n$ are called the {\em basis functions}. A very common choice of basis functions is the set of monomials $$\{g_0(t)=1,~g_i(t)=t^i,~i=1,\dots,n\},$$ so that the prototype is modelled as a polynomial function. In this study we are not limiting ourselves to polynomials, but require $\{g_i(t)\}_{i=0}^n$ to form a Chebyshev system,  that is
 \begin{equation}\label{eq:cheb_system}
 \det\{g_i(t_j)\},~i=0,\dots,n,~j=1,\dots,n,~a\leq t_1<\dots<t_{n+1}\leq b
 \end{equation} 
does not vanish for any choice of $t_k,~k=1,\dots,n$, such that $$a\leq t_1<\dots<t_{n+1}\leq b.$$

Our aim is to choose the set of parameters for the cluster prototype in such a way that the maximal deviation from each member of the group on $[a,b]$ is minimal.  That is, one has to solve the following optimisation problem:
\begin{equation}\label{eq:main_pol_convex}
{\rm minimise}~ F({\bf A})= \sup_{t_i\in [a,b],~i=1,\dots,N,~j=1,\dots,l}|S_j(t_i)-S({\bf A},t_i)|,
\end{equation} 
where ${\bf A}=(a_0,\dots,a_n)\in\R^{n+1}$, $a_k,~k=0,\dots,n$ are the approximation parameters and also the decision variables.
    
It can be shown that $F({\bf A})$ is convex, since it is a supremum of convex function. Therefore, we are working with an unconstrained convex problem with $n+1$ variables. It is also possible to formulate this problem as a linear programming problem. 

Let $$S_{{\rm max}}(t)=\max_{j=1,\dots,l}S_j(t)$$
and $$S_{{\rm min}}(t)=\min_{j=1,\dots,l}S_j(t).$$
Then (\ref{eq:main_pol_convex}) is equivalent to the following
$$
{\rm minimise}~ F({\bf A})= \sup_{t_i\in [a,b],~i=1,\dots,N}\max\{S_{\max}(t_i)-S({\bf A},t_i),S({\bf A},t_i)-S_{\min}(t_i)\}$$
and therefore the problem has been reduced to a {\em two curve  approximation problem}. 
Consider an additional variable $$z= \sup_{t_i\in [a,b],~i=1,\dots,N,~j=1,\dots,l}|S_j(t_i)-S({\bf A},t_i)|,$$ then the following linear programming problem is equivalent to~(\ref{eq:main_pol_convex}):
\begin{equation}\label{eq:main_pol_LP}
{\rm min}~z
\end{equation} 
subject to
\begin{align}\label{eq:constraints_pol_LP1}
S_{{\rm max}}(t_i)-S({\bf A},t_i)\leq z,~j=1,\dots,N;\\
\label{eq:constraints_pol_LP4}
S({\bf A},t_i)-S_{{\rm min}}(t_i)\leq z,~j=1,\dots,N.
\end{align} 
This linear programming problem has $n+2$ variables and $2N$ constraints. Since $N$ is the number of points where the signal segments are recorded, $N$ may be large. However, it is still more efficient to solve~(\ref{eq:main_pol_convex}) through its linear formulation.

 There are many efficient methods for solving linear programming problems. The first efficient linear programming algorithm (simplex method) was developed by G.~Dantzig in 1947~\cite{GB}. It was demonstrated in 1972 by~\cite{KM} that  the worst-case complexity of the simplex method is exponential. Despite this result, the simplex method is remarkably efficient and included in most linear programming packages.  Another important group of algorithms is interior point methods developed in~\cite{KN}, see also~\cite{JS}.  The  worst-case complexity of interior point methods is polynomial.
 
\subsection{Prototype update}
Suppose now that a signal group prototype $S^*(t)$ has been constructed and $$({\bf A}^*,z^*)=(a^*_0,\dots,a^*_{n},z^*)$$ is the corresponding optimal solution to~(\ref{eq:main_pol_LP})-(\ref{eq:constraints_pol_LP4}). Let $S(t)$ be an additional signal segment that needs to be included. How can we recompute the group prototype.

 One way to proceed is to update $S_{{\rm max}}(t)$ and $S_{{\rm min}}(t)$ and solve another linear programming problem. There are a number of more efficient ways to approach this problem. In particular, the cluster centre does not have to be computed from scratch at each iteration. 
 
  Before moving forward, let us underline three obvious, but yet very important properties.  
 \begin{enumerate}
 \item If after moving a number of signal segments in and out of a curve cluster the corresponding $S_{\min}$ and $S_{\max}$ remain unchanged, then the cluster centre remains the same. 
\item  If after moving a number of signal segments in and out of a curve cluster the corresponding $S_{\min}$ and $S_{\max}$ change, but the points of the maximal deviation remain the same, then the cluster centre remains the same. 
 \item If there exists a point $t\in[a,b]$, such that 
  \begin{equation}\label{eq:update1}
  S_{\max}(t)-S_{\min}(t)=2\sup_{t\in[a,b]}\max\{S^*(t)-S_{\min}(t),S_{\max}(t)-S^*(t)\}
  \end{equation}
  then the prototype $S^*(t)$ does not require any update, since the approximation can not be improved.
  \end{enumerate}


Since the prototype update has to be recomputed repeatedly, one of our objectives is to demonstrate how the cluster prototype obtained at the previous iteration can be reused for the next one. This can be done, since the updated constraint matrix contains several rows from the previous iteration. A Sherman-Morrison formula-based approach for such kind of linear programs has been proposed in~\cite{Sukhorukova2015}.    

One possible approach is to use the solution from the previous iteration as an initial point for the   next one. This approach should be exercise with care, since the final point from the previous iteration may be infeasible for the next iteration, since there are several segments moving in and out the group. In section~\ref{sec:modelling} we propose a more robust approach. This approach is based on the well-known  de la  Vall{\'e}e-Poussin's procedure~\cite{valleepoussin:1911}, originally developed for classical polynomial approximation.

\section{Optimality conditions}\label{sec:opt_cond}

In this section we develop the necessary and sufficient optimality conditions for two curve approximation, that are based on convex analysis and alternating sequence. Since the objective function is convex, we will apply convex analysis approaches from~\cite{Rockafellar70, Zalinescu2002}.  Before proceeding to two curve approximation, we provide classical results of Chebyshev approximation.

Chebyshev approximation theory is concerned with the approximation of a function \(f\), defined on a (continuous or discrete) domain \(\Omega\), by another function \(s\) taken from a family \(\cal{F}\) (for example, polynomials of degree $n$). At any point \(t\in \Omega\) the difference 
\[d(t)\triangleq s(t)-f(t)\] is called the \emph{deviation} at \(t\), and the \emph{maximal absolute deviation} is defined as \[\|s-f\|\triangleq \sup_{t\in \Omega} |s(t) - f(t)|.\]
The problem of best Chebyshev approximation is to find a function \(s^*\in \cal{F}\) minimising the maximal absolute deviation over \(\cal{F}\). Such a function \(s^*\) is called a \emph{best approximation} of \(f\).

The seminal result of approximation theory is  Chebyshev's alternation theorem~\cite{chebychev1854theorie}. Let \(P_n\) be the set of polynomials of degree at most \(n\) with real coefficients.

\begin{theorem} (Chebyshev alternation theorem, 1854) \label{thm:Chebyshev}
    A polynomial \(p^*\in P_{n}\) is a best approximation to a continuous function \(f\) on an interval \([a,b]\) if and only if there exist \(n+2\) points
    \( a\leq t_1 < \ldots < t_{n+2} \leq b\) and a number \(\sigma\in \{-1,1\}\) such that 
    \[ (-1)^i \sigma (f(t_i) - p^*(t_i)) = \|f-p^*\|, \forall i=1\ldots,n+2. \]
\end{theorem}
The sequence of points \((t_i)_{i=1,\ldots,n+2}\) is called an \emph{alternating sequence}.

Recall that in our study we are not restricted to polynomials. We require cluster prototypes to have a form
\begin{equation}\label{eq:approximation_form}
S({\bf A},t)=\sum_{i=1}^{n}a_ig_i(t),
\end{equation}
 where vector ${\bf A}=(a_0,\dots,a_n)^T\in \R^{n+1}$ is the vector of parameters (decision variables) and functions $g_i(t)$ are the basis functions (given). The only requirement for the basis functions is to form a Chebyshev system in $[a,b]$.

Now we proceed to the Chebyshev approximation based curve clustering.
Supposed that a cluster consists of $m$ signals ($S_1,\dots,S_m$), assigned with respect to the shortest distance to the cluster centres.  Now we need to recompute the cluster prototype.

First of all, we need to construct two curves:
\begin{itemize}
\item $S_{\max}(t)=\max_{i=1,\dots,m} S_i(t)$;
\item $S_{\min}(t)=\min_{i=1,\dots,m} S_i(t)$.
\end{itemize}
Then the parameters of the cluster prototype are the solution of the following optimisation problem:
$${\rm minimise}~ F(A)= \max_{t\in[a,b]}\{S_{\max}(t)-S({\bf A},t),S({\bf A},t)-S_{\min}(t)\}.$$
If the interval $[a,b]$ is discretised ($t_i\in[a,b],~i=1,\dots,N$), a solution can be obtain by solving a linear programming problem  (see section~\ref{sec:math_formulation}).

Let
$$\Delta= \max_{t\in[a,b]}\{S_{\max}(t)-S({\bf A},t),S({\bf A},t)-S_{\min}(t)\}.$$
\begin{definition}
A point $t_k$ where $$S_{\max}(t_k)-S({\bf A},t_k)=\Delta~{\rm or}~S({\bf A},t_k)-S_{\min}(t_k)=\Delta$$ is called a maximal deviation point. 
\end{definition}
\begin{definition}
A maximal deviation point~$t_k$, such that $$S_{\max}(t_k)-S({\bf A},t_k)=\Delta$$ is called a positive deviation alternating point, while a point $t_k$,   such that $$S(A,t_k)-S_{\min}(t_k)=\Delta$$ is called a negative deviation alternating point.
\end{definition}
\begin{theorem}  \label{thm:two_curves}
    An approximation $S({\bf A^*},t)$ is a best approximation to a pair of curves $S_{\max}$ and $S_{\min}$ on an interval \([a,b]\) if and only if     at least one of the following conditions holds:
    \begin{enumerate}
    \item there exists a time moment $t_k\in[a,b]$, such that 
    $$\Delta=S_{\max}(t_k)-S({\bf A^*},t_k)=S({\bf A^*},t_k)-S_{\min}(t_k);$$
    \item     there exist \(n+2\) points
    \( a\leq t_1 < \ldots < t_{n+2} \leq b\) and
    \begin{itemize}
    \item $\Delta=S_{\max}(t_k)-S({\bf A^*},t_k)=S({\bf A^*},t_{k+1})-S_{\min}(t_{k+1}),~k=1,\dots,n$ or
    \item $\Delta=S({\bf A^*},t_{k})-S_{\min}(t_{k})=S_{\max}(t_{k+1})-S({\bf A^*},t_{k+1}),~k=1,\dots,n$.
    \end{itemize}
    \end{enumerate}
\end{theorem}
{\bf Proof:}

Since the objective function is convex, its necessary and sufficient optimality condition is as follows:
\begin{equation}
0_{n+1}\in \partial F({\bf A^*}),
\end{equation}
where ${\bf A^*}$ is an optimal set approximation parameters. The subdifferential
\begin{equation}
\label{eq:subdifferential}
 \partial F({\bf A^*})=\Delta{\rm co}\left\{\left( 
 \begin{matrix}
 g_0(t_k^+)\\
 g_1(t_k^+)\\
 \vdots\\
 g_n(t_k^+)\\
 \end{matrix}
 \right), -
 \left( 
 \begin{matrix}
 g_0(t_k^-)\\
 g_1(t_k^-)\\
 \vdots\\
 g_n(t_k^-)\\
 \end{matrix}
 \right)
 \right\},~t_k^+\in T^+,~t_k^-\in T^-,
\end{equation}
where $T^+$ is the set of positive deviation alternating points and $T^-$ is the set of negative deviation alternating points. 
 This condition is equivalent to the existence of  a positive solution of the following linear system:
\begin{equation}\label{eq:sol_matrix}
{\bf M\Lambda}=0_{n+1},
\end{equation} 
where ${\bf M}$ is a matrix whose columns are the gradients at the maximal deviation points (extreme points of the subdifferential~$ \partial F$) and ${\bf \Lambda}$ is a vector whose components are non-negative and the sum of all the components is~1 (that is, there is at least one strictly positive component). 

Due to Caratheodory\rq{}s theorem, there exists a system of at most $n+2$  points from the subdifferential whose convex combination gives~$0_{n+1}$. Therefore, the number of columns in matrix ${\bf M}$ is at most $n+2$. Since the basis functions $g_i(t),~i=1,\dots,n$ form a Chebyshev system, the number of columns in ${\bf M}$ can not be less than $n+2$, otherwise the only solution to~(\ref{eq:sol_matrix}) is the trivial solution (all components are zeros). 

First assume that there is no point where both positive and negative maximal deviation is reached. Arrange maximal deviation points in ascending order:  $$a\leq t_1<t_2<\dots<t_{n+2}\leq b.$$
Hence, the system is as follows:
\begin{equation}\label{eq:system_optim}
\left(\sigma_0\left( 
 \begin{matrix}
 g_0(t_1)\\
 g_1(t_1)\\
 \vdots\\
 g_n(t_1)\\
 \end{matrix}
 \right)
\dots \sigma_n\left( 
 \begin{matrix}
 g_0(t_{n+1})\\
 g_1(t_{n+1})\\
 \vdots\\
 g_n(t_{n+1})\\
 \end{matrix}
 \right)\right)
 \left( 
 \begin{matrix}
 \alpha_0\\
 \alpha_1\\
 \vdots\\
 \alpha_n\\
 \end{matrix}
 \right)
=-\sigma_{n+1}\alpha_{n+1} \left( 
 \begin{matrix}
 g_0(t_{n+2})\\
 g_1(t_{n+2})\\
 \vdots\\
 g_n(t_{n+2})\\
 \end{matrix}
 \right),
\end{equation}
where $\sigma_{k}=1$ for positive alternating points and $\sigma_{k}=-1$ for negative alternating points ($k=0,\dots,n+1$). Since the functions $g_i(t),~i=0,\dots,n$ form a Chebyshev system (the corresponding determinants do not vanish and therefore the sign remains unchanged), applying Cramer\rq{}s rules, obtain that $$\sigma_i=-\sigma_{i+1},~i=0,\dots,n-1.$$  

Now assume that there are points where both positive and negative maximal deviation are reached, that is  the first condition holds, the approximation can not be improved, since the maximal deviation can not be made any smaller than $$\Delta={1\over 2}(S_{\max}(t_k)-S_{\min}(t_k)).$$


This proves the theorem.

\hskip 300pt $\square$
  
Therefore, the results for classical Chebyshev approximation and uniform approximation based clustering are very similar to each other. In the next section we demonstrate that, despite all these similarities, there are several fundamental differences. Therefore, a careful analysis is required for the extension of the classical results to the case of uniform approximation based clustering.
 
\section{Modelling}\label{sec:modelling}

In the case of classical Chebyshev approximation there are two important properties.
\begin{enumerate}
\item If the basis functions  form a Chebyshev system, then the optimal solution is unique.
\item By increasing the degree of the polynomial the maximal error can be made arbitrary small.
\end{enumerate}

The following simple example demonstrate that these two properties are not true for the case of  uniform approximation based clustering (two simultaneous curves approximation).

\begin{example}
Let $[a,b]=[0,1]$, $S_{\max}(t)=1-0.5t$ and $S_{\min}(t)=0.5t$. Find a best linear approximation for these two curves.

$S(A,t)=0.5$ is optimal. Moreover, any line
$$0.5+kt,~k\in[-0.25,0.25]$$ 
is optimal.

Regardless of the degree of the polynomial, the maximal deviation can not be made below 0.5. 
\end{example}

The highlighted fundamental differences between the classical uniform approximation and two curves simultaneous approximation demonstrate that not all the classical results can be generalised.   In the rest of this section we show that one of the fundamental results of the classical Chebyshev approximation, namely, de la  Vall{\'e}e-Poussin's procedure can be extended.

The classical de la  Vall{\'e}e-Poussin's procedure for polynomial approximation~\cite{valleepoussin:1911} has been extended to any basis functions, providing that they form a Chebyshev system (see~\cite{Karlin66} for details).  One starts with an initial basis (a system of $n+2$ points from $[a,b]$) then updates this basis by replacing one or more basis points by some other points from $[a,b]$ and eventually constructs an approximation that satisfies the necessary and sufficient optimality conditions. This procedure consists of two main steps.
 \begin{enumerate}
 \item Construct an approximation that deviates at the basis points from the original function by the same absolute deviation, the signs of the deviations are alternating. This approximation is also called the {\em Chebyshev interpolation approximation} (also Chebyshev interpolation polynomial). In the case of Chebyshev systems, such an approximation is unique~\cite{Karlin66, remez57}. 
 \item If there exists a non-basis point whose maximal absolute deviation is greater than it is at the basis points, then this point should be included into the basis, while one of the basis points should be removed. Namely, the removal  has to be done in such a way that the deviation signs at the new basis are alternating. It is enough to remove a neighbouring basis point with the same deviation sign or, if there is only one neighbouring point and the deviation sign at this point is opposite, remove the most extreme basis point from the opposite side of the interval.  
 \end{enumerate} 
By repeating these steps, one eventually arrives to the situation where there is no point that should be moved into the basis and therefore the current approximation is optimal, since the necessary and sufficient optimality conditions are satisfied. The basis exchange rule implies that each updated basis leads to the Chebyshev interpolation approximation, whose absolute deviation at the basis points is at least as large as the absolute deviation at the previous basis.

In the case of two curve approximation, the definition of basis remains the same: any set of distinct $n+2$ points from $[a,b]$ forms a basis.
Basis points  are also called {\em nodes}.

In our study we are working with two curves and therefore there may be points whose absolute deviation is maximal and both deviation signs are active. To be able to work with this kind of points, we introduce the following definitions.
\begin{definition}
A node where the absolute deviation is maximal and both positive and negative deviation is reached is called a double node.
\end{definition}

\begin{lemma}\label{lem:trVP}
Assume that the following two properties hold.
\begin{enumerate}
\item An original basis is chosen in such a way that there is no double node.
\item The replacement of one of the basis point by  a maximal absolute deviation point leads to a basis without double nodes.
\end{enumerate}
Then the classical  de la  Vall{\'e}e-Poussin's procedure can be extended to the case of two curve approximation for Chebyshev systems and terminate at a point where the condition~2 of  Theorem~\ref{thm:two_curves} is satisfied.
\end{lemma}
{\bf Proof:}

We start with an arbitrary collection of $n+2$ points from~$[a,b]$ as the initial basis~$T^0=\{t^0_0,\dots,t^0_{n+1}\}$. For each point $t\in T^0$ from this basis assign one of the values $S_{\max}(t)$ or $S_{\min}(t)$ in such a way that there is no pair of  neighbouring points assigned to the same curve (that is, construct $F(t)$). In this case there exists a unique Chebyshev interpolation approximation that deviates at the basis points from the assigned values by the same absolute value and the signs of the deviations are alternating (similar to one curve approximation).

The generalisation of the basis update step requires to consider two different possibilities.
\begin{enumerate}
\item The absolute maximal deviation outside of the basis does not exceed the absolute deviation at the basis points. In this case the condition~2 of Theorem~\ref{thm:two_curves} is satisfied and therefore the obtained approximation is optimal.
\item There exists a point outside of the basis where the absolute deviation is higher than it is at the basis points. The basis update rule is the same as it is for one curve approximation:  the maximal absolute deviation point replaces the adjacent  basis point with the same deviation sign or, if there is no same sign adjacent basis point, this point replaces    the furthest basis point regardless of its deviation sign.  Similar to one curve approximation~\cite{Karlin66}, this basis update leads to a Chebyshev interpolation approximation with a higher absolute deviation. This can be demonstrated by using~$F(t)$.
%
\end{enumerate} 
%
%
%
%
Since we only consider the situation where there is no point which is both positive and negative alternation point, the procedure is fully extended.

\hskip300pt $\square$

Our next step is to demonstrate how the assumptions from Lemma~\ref{lem:trVP} can be removed. In Theorem~\ref{thm:VP} we demonstrate how to construct an optimal approximation where the condition~2 of Theorem~\ref{thm:two_curves} is satisfied or, if this is not possible,  how to find a basis, such that the assumptions of Lemma~\ref{lem:trVP} are satisfied. 

\begin{theorem}\label{thm:VP}
 The classical  de la  Vall{\'e}e-Poussin's procedure can be extended to the case of two curve approximation for Chebyshev systems.
\end{theorem}
{\bf Proof:}
We start by identifying points $t$, where the difference $$S_{\max}(t)-S_{\min}(t)$$ reaches its maximal value in~$[a,b]$. We will call these points maximal difference points. Note that there exists at least one maximal difference point. Let
$$\Delta^*=0.5\max_{t\in[a,b]}(S_{\max}(t)-S_{\min}(t)),$$
where $\Delta^*$ is a lower bound for the optimal maximal deviation.

Consider two possibilities.
\begin{enumerate}
\item   Assume that the number of maximal difference points is~$l\geq n+1$. Choose any $n+1$ maximal distance points.  Since the basis functions form a Chebyshev system, there exists a unique approximation that passes through 
$$0.5(S_{\max}(t)-S_{\min}(t))$$
at the chosen $n+1$ maximal distance points. If there is no point whose maximal absolute deviation is strictly higher than~$\Delta^*$, then the obtained approximation is optimal. 

Now assume that there exists a point~$t^*$ whose absolute maximal deviation is greater than~$\Delta^*$. This point, together with the chosen $n+1$~points, form a basis and the alternation order between~$S_{\max}(t)$ and $S_{\min}(t)$ is determined by the maximal absolute deviation at~$t^*$. By construction, this point can not have both positive and negative maximal deviation and therefore the alternation order is determined uniquely.  The absolute maximal deviation is higher than~$\Delta^*$, therefore, there is no double node and any possible basis update can not lead to the presence of double nodes. By Lemma~\ref{lem:trVP}, the procedure is extended. 

\item Assume that the number of maximal difference points~$l\leq n$. Add any arbitrary chosen $n+1-l$ points. Assign all these $n+1$ points to $S_{\max}(t)$ or $S_{\min}(t)$ in an alternating way. Since the basis functions form a Chebyshev system, there exists a unique approximation that passes through the specified $(n+1)$~points, the absolute deviation is $\Delta^*$. If there is no points where the maximal absolute deviation exceeds~$\Delta^*$, the current approximation is optimal. Otherwise, the maximal deviation point together with the current $(n+1)$~points form a basis ($n+2$ points in total, since none of the points can coincide). Therefore, there exists a unique Chebyshev interpolation approximation, whose absolute deviation exceeds~$\Delta^*$ at the basis points.   Therefore, there is no double node and any possible basis update can not lead to the presence of double nodes. By Lemma~\ref{lem:trVP}, the procedure is extended. 
\end{enumerate}
Therefore, the procedure is fully extended.

\hskip300pt  $\square$

Coming back to the cluster prototype update, the basis point, obtained at the previous iteration can be used as an initial basis for the following one. This approach is similar to the usage of the solution obtain at a previous iteration as an initial point for the next one, but there is no risk of getting an infeasible point. The proposed approach is very efficient when the number of curves moving in or out the curve cluster is not very large.

\section{Conclusions and further research directions}\label{sec:conc}

This paper extends the classical Chebyshev approximation results to the case of curve clustering. The idea to extend these results comes from the application of optimisation and approximation to signal clustering and can be viewed as an extension of the well-known $k$-means method to curve clustering in uniform metric ($k$-medoid). Therefore, apart from the theoretical significance for optimisation and approximation, the results have potential applications in the area of signal processing and other areas of engineering and science.

There are still a number of open problems. The following  conjecture is one of our future research directions.
\begin{conjecture}
The two curve approximation de la  Vall{\'e}e-Poussin's procedure is equivalent to  the dual simplex method applied to the corresponding linear programming problem.
\end{conjecture}
This conjecture is of both theoretical and applied significance, since it underlines an efficient approach for curve clustering.

 We also intend to study the theoretical properties of the original objective function for $k$-medoids, similar the study proposed in~\cite{BagUgonCluster} for the Euclidean norm-based distances.





    \bibliographystyle{amsplain}

\end{document}